\documentclass[11pt,english]{article}
\usepackage[T1]{fontenc}
\usepackage[latin1]{inputenc}
\usepackage{amsmath}
\usepackage{amssymb}

\makeatletter

\title{The Semigroup of a Word}
\author{Peter M. Higgins \& Norman R. Reilly}
\date{}
\def\qed{\quad\vrule height4.17pt width4.17pt depth0pt}

\makeatother

\usepackage{babel}
\begin{document}

\title{\textbf{Involution matchings, the semigroup of orientation-preserving
and orientation-reversing mappings, and inverse covers of the full
transformation semigroup}}

\author{Peter M. Higgins, University of Essex, U.K.}
\maketitle
\begin{abstract}
We continue the study of permutations of a finite regular semigroup
that map each element to one of its inverses, providing a complete
description in the case of semigroups whose idempotent generated subsemigroup
is a union of groups. We show, in two ways, how to construct an involution
matching on the semigroup of all transformations which either preserve
or reverse orientation of a finite cycle. Finally, by way of application,
we prove that when the base set has more than three members, a finite
full transformation semigroup has no cover by inverse subsemigroups
which is closed under intersection. 
\end{abstract}

\section{Introduction and General Results}

\subsection{Background}

In {[}6{]} the author introduced the study of \emph{permutation matchings,
}which are permutations on a finite regular semigroup $S$ that map
each element to one of its inverses. It follows from Hall's Marriage
Lemma that $S$ will possess a permutation matching if and only if
$S$ satisfies the condition that $|A|\leq|V(A)|$ for all subsets
$A$ of $S$ with set of inverses $V(A)$. Although not all finite
regular semigroups have a permutation matching, there are positive
results for many important classes. In {[}7{]} the author characterised
some classes of finite regular semigroups by the nature of their permutation
matchings and determined, in terms of Green's relations on principal
factors, when a finite orthodox semigroup $S$ has a permutation matching.
In this case a permutation matching implies the existence of an involution
matching. In Section 1.3 we show how this result may be extended to
semigroups whose idempotent-generated subsemigroup is a union of groups. 

It is not known whether the semigroup ${\cal O}_{n}$ of all order-preserving
mappings on a finite $n$-chain has a permutation matching of any
kind. It was shown in {[}6{]} however that ${\cal OP}_{n}$, the semigroup
of all orientation-preserving mappings on an $n$-cycle, has a natural
involution matching. In Section 2.1 we summarise relevant properties
of this semigroup and of ${\cal P}_{n}$, the semigroup of all orientation-preserving
and orientation-reversing mappings on an $n$-cycle. This latter semigroup,
which was introduced in {[}1{]} and independently by McAlister in
{[}9{]}, has an intricate structure, which is manifested in the context
of the problem of this paper. In Section 3 we construct a dual pair
of involution matchings of ${\cal P}_{n}$. 

There are no known examples of a finite regular semigroup $S$ that
has a permutation matching but no involution matching. It was proved
in {[}6{]} by graph theoretic techniques that ${\cal T}_{n}$, the
full transformation semigroup on an $n$-set, has a permutation matching
but it is not known if ${\cal T}_{n}$ has an involution matching.
However in Section 4 we show that ${\cal T}_{n}$ $(n\geq4)$ has
no involution matching through so-called strong inverses, which allows
us to show that ${\cal T}_{n}$ $(n\geq4)$ has no cover by inverse
semigroups that is closed under intersection. 

Following the texts {[}8{]} and {[}5{]}, we denote the set of idempotents
of a semigroup $S$ by $E(S)$. We shall write $(a,b)\in V(S)$ if
$a$ and $b$ are mutual inverses in $S$ and denote this as $b\in V(a)$
so that $V(a)$ is the set of inverses of $a\in S$. We extend the
notation for inverses to sets $A$: $V(A)=\bigcup_{a\in A}V(a)$.
Standard results on Green's relations, particularly those stemming
from Green's Lemma, will be assumed (Chapter 2 of {[}5{]}, specifically
Lemma 2.2.1) and fundamental facts and definitions concerning semigroups
that are taken for granted in what follows are all to be found in
{[}5{]}. We shall sometimes write ${\cal G}$ to stand for either
of the Green's relations ${\cal L}$ or ${\cal R}$. 

We say that a semigroup $S$ is \emph{combinatorial }(or \emph{aperiodic})
if Green's ${\cal H}$-relation on $S$ is trivial. A completely $0$-simple
combinatorial semigroup is known as a \emph{$0$-rectangular band}.
The full transformation semigroup on a base set $X$ is denoted by
${\cal T}_{X}$ or by ${\cal T}_{n}$ when $X=X_{n}=\{0,1,2,\cdots,n-1\}$. 

Let $C=\{A_{i}\}_{i\in I}$ be any finite family of finite sets (perhaps
with repetition of sets). A set $\tau\subseteq\bigcup A_{i}$ is a
\emph{transversal }of $C$ if there exists a bijection $\phi:\tau\rightarrow C$
such that $t\in\phi(t)$ for all $t\in\tau$. We assume Hall's Marriage
Lemma in the form that $C$ has a transversal if and only if \emph{Hall's
Condition }is satisfied, which says that for all $1\leq k\leq|I|$,
the union of any $k$ sets from $C$ has at least $k$ members. 

\subsection{Permutation matchings}

\textbf{Definitions 1.2.1 }Let $S$ be any semigroup and let $F=\{f\in T_{S}:f(a)\in V(a)\,\forall a\in\mbox{dom\,\ \ensuremath{f}\}}.$
We call $F$ the set of \emph{inverse matchings }of $S$. We call
$f\in F$ a \emph{permutation matching }if $f$ is a permutation of
$S$; more particularly $f$ is an \emph{involution matching} if $f^{2}=\varepsilon$,
the identity mapping on $S$. 

In the remainder of the paper we shall assume that $S$ is regular
and finite unless otherwise indicated. We shall often denote a matching
simply by $'$, so that the image of $a$ is $a'$. We use the shorthand
$a''$ as an abbreviation for $(a')'$. We shall work with the family
of subsets of $S$ given by $V=\{V(a)\}_{a\in S}$. The members of
$V$ may have repeated elements\textemdash for example $S$ is a rectangular
band if and only if $V(a)=S$ for all $a\in S$. However, we consider
the members of $V$ to be marked by the letter $a$, so that $V(a)$
is an unambiguous member of $V$ (strictly, we are using the pairs
$\{a,V(a)\},$ $(a\in S)$). We summarise some results of {[}5{]}.

\textbf{Theorem 1.2.2 }{[}6{]} For a finite regular semigroup $S$
the following are equivalent:

(i) $S$ has a permutation matching;

(ii) $S$ is a transversal of $V=\{V(a)\}_{a\in A}$; 

(iii) $|A|\leq|V(A)|$ for all $A\subseteq S$; 

(iv) $S$ has a permutation matching that preserves the ${\cal H}$-relation;
(meaning that $a{\cal H}b\Rightarrow a'{\cal H}b'$);

(v) each principal factor $D_{a}\cup\{0\}$ $(a\in S)$ has a permutation
matching; 

(vi) each $0$-rectangular band $B=(D_{a}\cup\{0\})/{\cal H}$ has
a permutation matching.

In {[}7, Remark 1.5{]} it was shown that we may replace `permutation
matching' by `involution matching' in Theorem 1.2.2 as regards the
implications ((i) $\Leftrightarrow$ (v)) $\Leftarrow$ ((iv) $\Leftrightarrow$
(vi)) although the missing forward implication has not been resolved. 

\subsection{Permutation matchings for an E-solid semigroup}

\textbf{Definition 1.3.1 }A regular semigroup $S$ is defined to be
$E$-\emph{solid} if $S$ satisfies the condition that for all idempotents
$e,f,g\in E(S)$ 
\[
e{\cal L}f{\cal R}g\rightarrow\exists h\in E(S):e{\cal R}h{\cal L}f.
\]
An alternative characterisation of an $E$-solid semigroup is that
of a regular semigroup $S$ for which the idempotent-generated subsemigroup
$\langle E(S)\rangle$ is a union of groups {[}3, Theorem 3{]}.

We prove our result on $E$-solid semigroups via the corresponding
result for orthodox semigroups. The proof of this latter result involved
reducing the general problem to the case of $0$-rectangular bands
and then showing that the corresponding ${\cal D}$-class may be diagonalised
in that the ${\cal R}$- and ${\cal L}$-classes may be ordered so
that all idempotents are contained in rectangular blocks (which then
form the maximal rectangular subbands of $D$); $S$ then has a permutation
matching if and only if, within each ${\cal D}$-class of $S$, these
blocks are similar in the following sense. 

\textbf{Definition 1.3.2 }Let $U_{1}$ and $U_{2}$ be finite rectangular
bands, let $m_{i}$ and $n_{i}$ denote the respective number of ${\cal R}$-classes
and ${\cal L}$-classes of $U_{i}$ $(i=1,2)$. We say that $U_{1}$
and $U_{2}$ are \emph{similar }if $\frac{m_{1}}{n_{1}}=\frac{m_{2}}{n_{2}}$. 

\textbf{Theorem 1.3.3 }{[}7, Theorems 3.7 and 3.1{]} Let $S$ be a
finite orthodox semigroup. Then $S$ has a permutation matching if
and only if for each $0$-rectangular band $B=(D_{a}\cup\{0\})/{\cal H}$
$(a\in S)$ the maximal rectangular subbands of $B$ are pairwise
similar. In that case the permutation matching of $S$ may be chosen
to be an involution matching.

\textbf{Proposition 1.3.4 }Each $0$-rectangular band $B=(D_{a}\cup\{0\})/{\cal H}$
$(a\in S)$ of a finite $E$-solid semigroup $S$ is orthodox. 

For completeness, we record a proof of the Proposition but note that
the class of all (not necessarily finite) $E$-solid semigroups is
a so-called \emph{e}-\emph{variety}, meaning that the class is\emph{
}closed under the taking of homomorphic images, of direct products,
and regular subsemigroups {[}4{]}. Also in {[}4{]} is shown that a
semigroup is orthodox if and only if the same is true of each of its
principal factors: (also see {[}5, Ex. 1.4.13(iv){]}).

\textbf{Proof} From the definition of $E$-solidity we see that each
principal factor $D_{a}\cup\{0\}$ of $S$ is itself $E$-solid, and
$B$ certainly is regular. Next we note that $B$ is $E$-solid through
two observations: $H\in E(B)$ if and only if $H=H_{e}$ for some
$e\in E(S)$, and $H_{a}{\cal G}H_{b}$ in $B$ if and only if $a{\cal G}b$
in $S$. Hence if $B$ contain three idempotents $H_{e},H_{f},$ and
$H_{g}$ with $e,f,g\in E(S)$, and they are such that $H_{e}{\cal L}H_{f}{\cal R}H_{g}$
in $B$, then $e{\cal L}f{\cal R}g$ in $S$ and by the $E$-solid
condition on $S$ we have $e{\cal R}h{\cal L}g$ for some $h\in E(S)$.
We now have $H_{e}{\cal R}H_{h}{\cal L}H_{f}$ in $B$ and $H_{h}\in E(B)$.
Therefore $B$ is $E$-solid.

To show that $B$ is indeed orthodox, first note that by our second
observation, $B$ is combinatorial. Then for any two idempotents of
$B$, which we now write as $e,f$, we have either that $ef=0\in E(B)$,
or otherwise $e{\cal L}g{\cal R}f$ for some $g\in E(B)$ whence,
since $B$ is $E$-solid, it follows that $h=ef\in E(B)$. Therefore
$E^{2}(B)=E(B)$, as required. $\qed$

\textbf{Theorem 1.3.5 }Let $S$ be a finite $E$-solid semigroup.
Then $S$ has a permutation matching if and only if the maximal rectangular
subbands of each of the $0$-rectangular bands $(D_{a}\cup\{0\})/{\cal H}$
are pairwise similar. Moreover if $S$ has a permutation matching
then $S$ has an involution matching. 

\textbf{Proof} By Theorem 1.2.2, $S$ has a permutation matching if
and only if the same can be said for all $B=(D_{a}\cup\{0\})/{\cal H}$
$(a\in S)$. By Proposition 1.3.4, each such $B$ is a finite orthodox
$0$-rectangular band. By Theorem 1.3.3, each such $B$ then has a
permutation matching if and only if the maximal rectangular subbands
of $B$ are pairwise similar, giving the first statement of Theorem
1.3.5. In this case, again by Theorem 1.3.3, each permutation matching
of each $B$ may be chosen to be an involution matching of $B$. Then
by (vi) implies (i) in Theorem 1.2.2 as it applies to involutions,
we conclude that $S$ itself has an involution matching, thus completing
the proof. $\qed$

\section{Matchings for ${\cal OP}_{n}$ and ${\cal P}_{n}$}

\subsection{The semigroups ${\cal OP}_{n}$ and ${\cal P}_{n}$}

We recap some of the important properties of the semigroups ${\cal OP}_{n}$
and ${\cal P}_{n}$. We also augment these results in order to build
a type of calculus for these semigroups. All semigroups under consideration
will be subsemigroups of ${\cal T}_{n}$. Basic properties of the
representation of $\alpha\in{\cal T}_{n}$ as a digraph $G(\alpha)$
can be found in the text {[}5, Section 1.5{]}. Each component $C$
of $G(\alpha)$ is \emph{functional}, meaning that each vertex has
out-degree $1$ so in consequence $C$ consists of a unique cycle
$Z(\alpha)$ with a number of directed trees rooted around the vertices
of $Z(\alpha)$. The set of cycle points of $G(\alpha)$ are exactly
the points in the \emph{stable range }of $\alpha$, denoted by stran$(\alpha)$,
which are the points of $X_{n}$ contained in the range of all powers
of $\alpha$. Pictures of these digraphs are helpful in seeing what
is going on and the reader is invited to draw them where relevant,
especially in the examples of Section 4 where they are a natural aid
to understanding.

For $\alpha\in{\cal T}_{n}$ we write $R=R(\alpha)$ for its range
$X\alpha$, while $t=|R(\alpha)|$ will stand for the \emph{rank}
of $\alpha$. The kernel relation of $\alpha$ on $X$ will be denoted
as ker$(\alpha)$ with the corresponding partition of $X_{n}$ written
as Ker$(\alpha)$. The set of fixed points of $\alpha$ will be denoted
by $F(\alpha)$. Facts from the source paper {[}1{]} are listed using
the term Result.

\textbf{Definitions 2.1.1 }(i) the \emph{cyclic interval }$[i,i+t]$
$(0\leq t\leq n-1)$ is the set $\{i,i+1,\cdots,i+t\}$ if $i+t\leq n-1$
and otherwise is the set \newline $\{i,i+1,\cdots,n-1,0,1,\cdots,(i+t)\,\mbox{(mod\,\ensuremath{n)\}}}$.

(ii) A finite sequence $A=(a_{0},a_{1},\cdots,a_{t})$ from $[n]$
is \emph{cyclic }if there exists no more than one subscript $i$ such
that $a_{i}>a_{i+1}$ (taking $t+1=0$). We say that $A$ is \emph{anti-cyclic
}if the reverse sequence $A^{r}=(a_{t},a_{t-1},\cdots,a_{0})$ is
cyclic. 

\textbf{Remarks 2.1.2 }To say that $A$ is cyclic as in (ii) is equivalent
to saying that for some subscript $i$, $a_{i+1}\leq\cdots\leq a_{t}\leq a_{0}\leq\cdots\leq a_{i}$
and the subscript $i$ with this property is unique unless $A$ is
constant. On the other hand $A$ is anti-cyclic means $A^{r}$ is
cyclic so that $A$ is anti-cyclic if and only if for some subscript
$i$ we have $a_{i+1}\geq\cdots\geq a_{t}\geq a_{0}\geq\cdots\geq a_{i}$
(and $i$ is unique if $A$ is not constant). The properties of cyclicity
and anti-cyclicity are inherited by subsequences and by sequences
obtained by cyclic re-ordering.

\textbf{Definition 2.1.3 }A mapping $\alpha\in{\cal T}_{n}$ is \emph{orientation-preserving
}if its list of images, $(0\alpha,1\alpha,\cdots,(n-1)\alpha)$, is
cyclic. The collection of all such mappings is denoted by ${\cal OP}_{n}$.\textbf{
}We say that $\alpha\in{\cal T}_{n}$ is \emph{orientation-reversing
}if $(0\alpha,1\alpha,\cdots,(n-1)\alpha)$ is anti-cyclic and the
collection of all orientation-reversing mappings is denoted by ${\cal OR}_{n}$.

\textbf{Result 2.1.4} ${\cal OP}_{n}$ is a regular submonoid of ${\cal T}_{n}$.
Each kernel class of $\alpha\in{\cal OP}_{n}$ is a cyclic interval
of $[n]$ and the maximal cycles of the components of the digraph
$G(\alpha)$ have the same number of vertices, denoted by $c(\alpha)$. 

\textbf{Definition 2.1.5} Let $\alpha\in{\cal OP}_{n}$ be of rank
$t\geq2$. We index the members of Ker$(\alpha)$ as $K_{i}$ $(0\leq i\leq t-1)$
in such a way that the set of initial points $a_{i}$ of the cyclic
intervals $K_{i}$ satisfy $a_{0}<a_{1}<\cdots<a_{t-1}$, denoting
this ordered set by $K(\alpha)$. The list $\{K_{0},K_{1},\cdots,K_{t-1}\}$
is called the \emph{canonical listing} of the kernel classes of $\alpha$.
For $r_{i}\in R(\alpha)$ where $r_{0}<r_{1}<\cdots<r_{t-1}$ we denote
the cyclic interval $[r_{i},r_{i}+1,\cdots,r_{i+1}-1]$ by $R_{i}$. 

\textbf{Result 2.1.6 }({[}1{]}, Theorem 3.3) For $t\geq2$ there is
a one-to-one correspondence $\Phi_{0}$ between the set of triples
$(K,R,i)$ where $K$ and $R$ are ordered $t$-sets of $X_{n}$ $(0\leq i\leq t-1)$
and $\{\alpha\in{\cal OP}_{n}:|X\alpha|=t\}$ whereby $(K,R,i)\mapsto\alpha$,
where each $a_{j}\in K$ is an initial point of a kernel class of
$\alpha$ and $a_{j}\alpha=r_{i+j}$$\,(0\leq i\leq t-1)$, subscripts
calculated modulo $t$. Moreover $H_{\alpha}=\{\Phi_{0}(K,R,i):\,i=0,1,\cdots,t-1\}$,
and so $|H_{\alpha}|=t$. 

\textbf{Result 2.1.7 }(i) The collection ${\cal P}_{n}={\cal OP}_{n}\cup{\cal OR}_{n}$
is a regular submonoid of ${\cal T}_{n}$; ${\cal R}$- ${\cal L}$-
and ${\cal D}$-classes are determined by equality or kernels, of
images, and of ranks respectively (as in ${\cal T}_{n}$ and ${\cal OP}_{n}$). 

(ii) the \emph{reflection mapping} $\gamma:[n]\rightarrow[n]$, whereby
$i\mapsto n-i-1$ $(i\in[n])$ is orientation-reversing and ${\cal P}_{n}=\langle a,e,\gamma\rangle$,
where $a$ is the $n$-cycle $(0\,1\,\cdots\,n-1)$ and $e$ is any
idempotent in ${\cal OP}_{n}$ of rank $n-1$; $\langle a,e\rangle={\cal OP}_{n}$.

(iii) ${\cal OP}_{n}\cap{\cal OR}_{n}=\{\alpha\in{\cal OP}_{n}:\,\mbox{rank\ensuremath{(\alpha)\leq2\}}. }$

(iv) $({\cal OR}_{n})^{2}={\cal OP}_{n}$, ${\cal OP}_{n}\cdot{\cal OR}_{n}={\cal OR}_{n}\cdot{\cal OP}_{n}={\cal OR}_{n}$.

It is also proved in {[}1{]} and in {[}9{]} that the respective maximal
subgroups of rank $t$ of ${\cal OP}_{n}$ and of ${\cal P}_{n}$
are cyclic groups of order $t$ and dihedral groups of order $2t$.
Also every non-constant member $\alpha\in{\cal OP}_{n}$ factorizes
uniquely as $\alpha=a^{r}\phi$ where $a$ is the $n$-cycle as above
and $\phi\in{\cal O}_{n}$. The constant mappings on $[n]$ comprise
$D_{1}$, the lowest ${\cal D}$-class of ${\cal P}_{n}$. Any permutation
of $D_{1}$ is a permutation matching of $D_{1}$ and for that reason
$D_{1}$ will not need to feature in our subsequent discussion. 

\textbf{Definition 2.1.8 }For $\alpha\in{\cal P}_{n}$ we shall write
$\rho(\alpha)=(K,R)$, where $K$ and $R$ are the respective sets
$K(\alpha)$ of initial points of kernel classes and $R(\alpha)$. 

Note that for any $\alpha,\beta\in S={\cal OP}_{n}$ or ${\cal P}_{n}$,
$\alpha{\cal H}\beta$ if and only if $\rho(\alpha)=\rho(\beta)$.
We now extend Result 2.1.6 to ${\cal P}_{n}$. 

\textbf{Theorem 2.1.9 }For $t\geq2$ there is a one-to-one correspondence
$\Phi$ between the set of quadruples $(K,R,i,k)$ where $K$ and
$R$ are ordered $t$-sets of $[n]$, $0\leq i\leq t-1$, $k=\pm1$
and $\{\alpha\in{\cal P}_{n}:|X\alpha|=t\}$. The correspondence is
given by $(K,R,i,k)\mapsto\alpha$, where each $a_{j}\in K$ is an
initial point of a kernel class of $\alpha$ and $a_{j}\alpha=r_{i+kj}$$\,(0\leq j\leq t-1)$,
subscripts calculated modulo $t$. 

\textbf{Proof} The first statement of Result 2.1.6 implies that $\Phi|_{k=1}$
maps bijectively onto the set of non-constant mappings in ${\cal OP}_{n}$.
We show that $\Phi|_{k=-1}$ maps bijectively onto the set of non-constant
mappings in ${\cal OR}_{n}$. The equation $K_{j}\alpha=r_{i-j}$
certainly specifies a unique mapping $\alpha=\Phi(K,R,i,-1)\in{\cal T}_{n}$,
and distinct quadruples yield distinct mappings. We need to check
that $\alpha\in{\cal OR}_{n}$. We have however the following equality
of two lists: 
\begin{equation}
a_{i+1}\alpha=r_{t-1}>a_{i+2}\alpha=r_{t-2}>\cdots>a_{i}\alpha=r_{0}
\end{equation}
It follows from (1) and Remarks 2.1.2 that the image of the cyclic
list $K$ under $\alpha$ is anti-cyclic and so $\alpha=\Phi(K,R,i,-1)\in{\cal OR}_{n}$;
hence $\Phi|_{k=-1}$ is a one-to-one mapping into the set of mappings
of ${\cal OR}_{n}$ of rank at least $2$. 

Conversely, let $\alpha\in{\cal OR}_{n}$ be of rank $t\geq2$. Since
multiplication on the right by $\gamma$ defines a bijection of ${\cal OP}_{n}$
onto ${\cal OR}_{n}$, it follows that the kernel classes of $\alpha$
are cyclic intervals and so $H_{\alpha}$ is determined by a pair
of ordered $t$-sets $(K,R)$. Take $i$ such that $a_{i}\alpha=r_{0}$.
Then since $a_{i+1},a_{i+2},\cdots,a_{i}$ is cyclic and $\alpha\in{\cal OR}_{n}$,
it follows that $a_{i+1}\alpha,a_{i+2}\alpha,\cdots,a_{i}\alpha$
is anti-cyclic. However, since $a_{i}\alpha=r_{0}=$ min$R$, it follows
by Remarks 2.1.2 that (1) holds for $\alpha$ and so $\alpha=\Phi(K,R,i,-1)$.
Therefore $\Phi|_{|k=-1}$ is a bijection onto the set of non-constant
mappings of ${\cal OR}_{n}$. Finally note that for $k=\pm1$, the
rank of $\alpha=\Phi(K,R,i,k)$ is indeed $t=|R|=|K|$. $\qed$

\textbf{Corollary 2.1.10 }For $t\geq3$, each ${\cal H}$-class $H$
of ${\cal P}_{n}$ contained in $D(t)$ is a disjoint union $H=(H\cap{\cal OP}_{n})\cup(H\cap{\cal OR}_{n})$
with each set in the union of cardinal $t$. 

\textbf{Proof} Let $H=\{\alpha\in{\cal P}_{n}:\rho(\alpha)=(K,R)\}$.
Then by Theorem 2.1.9, $H\cap{\cal OP}_{n}=\{\Phi(K,R,i,1):0\leq i\leq t-1\}$
and $H\cap{\cal OR}_{n}=\{\Phi(K,R,i,-1):0\leq i\leq t-1\}$; these
two sets each have $t$ members and are disjoint by Result 2.1.7(iii).
$\qed$

We shall refer to the coding of each $\alpha\in{\cal P}_{n}$ in the
form $(K,R,i,k)$ as the \emph{KRik-coordinates }of $\alpha$, noting
that $(K,R,i,k)=\Phi^{-1}(\alpha)$. We call $i$ and $k$ respectively
the \emph{shift }and the \emph{parity }of $\alpha$. 

\textbf{Lemma 2.1.11 }Let $\alpha\in{\cal P}_{n}$ with $\rho(\alpha)=(K,R)$.
Then 

(i) $\rho(\alpha\gamma)=(K,n-1-R)$; ~~~(ii) $\rho(\gamma\alpha)=(n-K,R)$; 

(iii) $\rho(\gamma\alpha\gamma)=(n-K,n-1-R)$.

\textbf{Proof }(i) is immediate from definition as is the fact that
$R(\gamma\alpha)=R(\alpha)$ in (ii). Continuing in (ii), suppose
that $\alpha=\Phi(K,R,i,k)$. Then for $r_{i+kj}\in R$ $(0\leq j\leq t-1)$,
we obtain: 
\[
r_{i+kj}(\gamma\alpha)^{-1}=r_{i+kj}\alpha^{-1}\gamma^{-1}=K_{j}\gamma
\]
\[
=(n-1)-\{a_{j},a_{j}+1,\cdots,a_{j+1}-1\}=\{n-a_{j+1},n-a_{j+1}+1,\cdots,n-a_{j}-1\},
\]
where $j+1$ is calculated modulo $t$. It follows that $K(\gamma\alpha)=n-K(\alpha)$,
thereby establishing (ii). Applying (i) and then (ii) now gives (iii)
as follows:
\[
\rho(\gamma\alpha\cdot\gamma)=(K(\gamma\alpha),n-1-R(\gamma\alpha))=(n-K,n-1-R)\,\,\qed
\]

\newpage

\textbf{Proposition 2.1.12} Let $\alpha=\Phi(K,R,i,k)\in{\cal P}_{n}$.
Then

(i) $\alpha\gamma=\Phi(K,n-1-R,-(i+1),-k)$;

(ii) if $0\not\not\in K$ then $\gamma\alpha=\Phi(n-K,R,i-2k,-k)$; 

(iii) if $0\in K$ then $\gamma\alpha=\Phi(n-K,R,i-k,-k)$;

(iv) if $0\not\in K$ then $\gamma\alpha\gamma=\Phi(n-K,n-1-R,2k-(i+1),k);$

(v) if $0\in K$ then $\gamma\alpha\gamma=\Phi(n-K,n-1-R,k-(i+1),k)$. 

\textbf{Proof} (i) We are working throughout modulo $t$ on subscripts.
By Lemma 2.1.11(i) we have $\rho(\alpha\gamma)=(K,n-1-R)$. Now 
\[
n-1-R=\{n-1-r_{t-1}<n-1-r_{t-2}<\cdots<n-1-r_{0}\}.
\]
Let us denote $n-1-r_{-(j+1)}$ by $s_{j}$ $(0\leq j\leq t-1)$ so
that $R(\alpha\gamma)=\{s_{0}<s_{1}<\cdots<s_{t-1}\}$. Hence $a_{j}\alpha\gamma=r_{i+kj}\gamma=n-1-r_{i+kj}$;
now 
\[
i+kj=-(-i-kj),\,\text{so that\,\ensuremath{a_{j}\alpha\gamma}}=s_{-(i+1)-kj},
\]
which establishes equation (i). 

(ii) By Lemma 2.1.11(ii) we have $\rho(\gamma\alpha)=(n-K,R)$. Since
$1\leq a_{0}$ 
\[
n-K=n-a_{t-1}<n-a_{t-2}<\cdots<n-a_{0}.
\]
Let us denote $n-a_{-(j+1)}$ by $b_{j}$ $(0\leq j\leq t-1)$ so
that $K(\gamma\alpha)=\{b_{0}<b_{1}<\cdots<b_{t-1}\}$. Hence:
\[
b_{j}\gamma\alpha=(n-1-(n-a_{-(j+1)}))\alpha=(a_{-(j+1)}-1)\alpha=a_{-(j+2)}\alpha=r_{i-k(j+2)}=r_{(i-2k)-kj},
\]
which establishes equation (ii).

(iii) Now since $a_{0}=0$ we have $n-a_{0}=n\equiv0$ (mod $n)$
and so: 
\[
n-K=n-a_{0}<n-a_{t-1}<n-a_{t-2}<\cdots<n-a_{1}.
\]
Let us denote $n-a_{-j}$ by $b_{j}$ $(0\leq j\leq t-1)$ so that
$K(\gamma\alpha)=\{b_{0}<b_{1}<\cdots<b_{t-1}\}$. Hence 
\[
b_{j}\gamma\alpha=(n-1-(n-a_{t-j}))\alpha=(a_{-j}-1)\alpha=a_{-(j+1)}\alpha=r_{i-k(j+1)}=r_{(i-k)-kj},
\]
which establishes equation (iii).

(iv) By Lemma 2.1.11(iii) we have $\rho(\gamma\alpha\gamma)=(n-K,n-1-R)$.
Now using (ii) we obtain
\[
b_{j}\gamma\alpha\gamma=r_{(i-2k)-kj}\gamma=n-1-r_{-(2k-i+kj)}=s_{(2k-(i+1))+kj},
\]
which establishes equation (iv).

(v) By Lemma 2.1.11(iii) we have $\rho(\gamma\alpha\gamma)=(n-K,n-1-R)$.
Now using (iii) we obtain
\[
b_{j}\gamma\alpha\gamma=r_{i-k-kj}\gamma=n-1-r_{-(k-i+kj)}=s_{(k-(i+1))+kj},
\]
which establishes equation (v). $\qed$

\textbf{Example 2.1.13 }As an example we find $\gamma\alpha\gamma$
for $\alpha\in{\cal OR}_{10}$ given by:
\[
\alpha=\begin{pmatrix}0 & 1 & 2 & 3 & 4 & 5 & 6 & 7 & 8 & 9\\
3 & 2 & 2 & 8 & 8 & 6 & 6 & 4 & 3 & 3
\end{pmatrix};
\]
so that $n=10$, $t=5$, $K=\{1,3,5,7,8\}$, $R=\{2,3,4,6,8\}$ and
$\alpha=\Phi(K,R,0,-1)$. Since $0\not\in K$, according to Proposition
2.12(iv), we should find that $\gamma\alpha\gamma=\Phi(10-K,9-R,2,-1)$
as $i(\gamma\alpha\gamma)=2(-1)-(0+1)=-3\equiv2$ (mod $5$). Now
$n-K=\{2,3,5,7,9\}$, and $n-1-R=\{1,3,5,6,7\}$. This accords with
the direct calculation of $\gamma\alpha\gamma$, which corresponds
to reversing the image line of $\alpha$ (to get $\gamma\alpha)$
and subtracting the images from $n-1=9$.
\[
\gamma\alpha\gamma=\begin{pmatrix}0 & 1 & 2 & 3 & 4 & 5 & 6 & 7 & 8 & 9\\
6 & 6 & 5 & 3 & 3 & 1 & 1 & 7 & 7 & 6
\end{pmatrix}.
\]

\section{Permutation matchings for ${\cal OP}_{n}$ and ${\cal P}_{n}$ }

\subsection{An approach via subset involutions}

\textbf{Lemma 3.1.1 }Let $A\subseteq S$ and let $(\cdot')$ denote
an ${\cal H}$-class-preserving involution matching on the set $A$.
Then $(\cdot')$ may be uniquely extended to an involution matching
on $A_{{\cal H}}=\cup_{a\in A}H_{a}$. In particular, if $A$ meets
every ${\cal H}$-class of $S$ then $(\cdot')$ extends uniquely
to an ${\cal H}$-class-preserving involution of $S$, which we shall
call the \emph{induced involution matching }on $S$. 

\textbf{Proof }Since $(\cdot')$ is ${\cal H}$-class-preserving it
induces an involution on $A_{{\cal H}}$ by $H_{a}\mapsto H_{a'}$.
We then have an involution matching on $A_{{\cal H}}$ defined by
$b\mapsto b'$ where $b\in H_{a}$ say and $b'$ is the unique inverse
of $b$ in $H_{a'}$. $\qed$

We now construct what we shall refer to as the \emph{natural matching
involution }for ${\cal P}_{n}$, which is the induced involution matching
on ${\cal P}_{n}$ extending the involution of ${\cal OP}_{n}$ recorded
in {[}6{]}. 

\textbf{Theorem 3.1.2} The semigroup ${\cal P}_{n}$ has an ${\cal H}$-class-preserving
involution matching $(\cdot')$ defined by the rule:
\begin{equation}
(\alpha=\Phi(K,R,i,k))\Rightarrow(\alpha'=\Phi(R,K,-ki,k))
\end{equation}
\textbf{Proof }From Theorem 2.1.9 and its proof it follows that (2)
defines an ${\cal H}$-class-preserving involution that maps each
of ${\cal OP}_{n}$ and ${\cal OR}_{n}$ onto ${\cal OP}_{n}$ and
${\cal OR}_{n}$ respectively: certainly $\alpha''=\alpha$ as $(-k)(-ki)=k^{2}i=i.$
It remains only to check that $(\alpha,\alpha')\in V({\cal P}_{n})$
and by symmetry it is enough to verify that $\alpha=\alpha\alpha'\alpha$.
To this end take $x\in[n]$ with $x\in K_{j}$ say where the kernel
classes of $\alpha$ are labelled by subscripts in the canonical order.
Then since $-ki+k(i+kj)=-ki+ki+k^{2}j=j$ we obtain:
\[
x\alpha\alpha'\alpha=a_{j}\alpha\alpha'\alpha=r_{i+kj}\alpha'\alpha=a_{-ki+k(i+kj)}\alpha=a_{j}\alpha=x\alpha,
\]
and so $\alpha=\alpha\alpha'\alpha$, as required. $\qed$

Note that for $t=2$ we have ${\cal OP}_{n}\cap D_{2}={\cal OR}_{n}\cap D_{2}=D_{2}$
and in this case $k=-k$ (mod $t)$. This collapse in the fourth entry
of the \emph{KRik-}coordinates leads to the involution matching taking
on the simpler form $\alpha=\Phi(K,R,i)\mapsto\alpha'=\Phi(R,K,i)$. 

\textbf{Examples 3.1.3 }Let $n=8$, $t=4$, $K=\{0,2,4,6\}$, $R=\{1,3,5,7\}$
and $\alpha\in{\cal OR}_{8}$ defined by $\alpha=\Phi(K,R,3,-1)$
so that $\alpha'=\Phi(R,K,3,-1)$. 
\[
\alpha=\begin{pmatrix}0 & 1 & 2 & 3 & 4 & 5 & 6 & 7\\
7 & 7 & 5 & 5 & 3 & 3 & 1 & 1
\end{pmatrix}\,\,\alpha'=\begin{pmatrix}0 & 1 & 2 & 3 & 4 & 5 & 6 & 7\\
0 & 6 & 6 & 4 & 4 & 2 & 2 & 0
\end{pmatrix}.
\]
For an example in ${\cal OP}_{n}$ let us take $n=10,$ $t=6$, $K=\{0,2,4,7,8,9\},$
$R=\{0,1,2,5,6,7\}$ and $\text{\ensuremath{\alpha}}=\Phi(K,R,4,1)$.
Since $-i=-4=2$ (mod $6)$ we obtain $\alpha'=\Phi(R,K,2,1)$:
\[
\alpha=\begin{pmatrix}0 & 1 & 2 & 3 & 4 & 5 & 6 & 7 & 8 & 9\\
6 & 6 & 7 & 7 & 0 & 0 & 0 & 1 & 2 & 5
\end{pmatrix}\,\,\alpha'=\begin{pmatrix}0 & 1 & 2 & 3 & 4 & 5 & 6 & 7 & 8 & 9\\
4 & 7 & 8 & 8 & 8 & 9 & 0 & 2 & 2 & 2
\end{pmatrix}.
\]

\textbf{Example 3.1.4 }The natural involution inverse of a group element
is not necessarily a group element and nor does the natural involution
map ${\cal O}_{n}$ into itself. Both features are seen in the following
example. Take $n=3$, $t=2$, $K=\{0,2\}$, $R=\{1,2\}$, and put
$\alpha=\Phi(K,R,0,1)$ so that $\alpha'=\Phi(R,K,0,1)$: 
\[
\alpha=\begin{pmatrix}0 & 1 & 2\\
1 & 1 & 2
\end{pmatrix}\,\,\alpha'=\begin{pmatrix}0 & 1 & 2\\
2 & 0 & 2
\end{pmatrix};
\]
we see that $\alpha\in E({\cal O}_{3})\subseteq E({\cal OP}_{3})$
and so is an order-preserving group element while $\alpha'^{2}$ is
the constant mapping with range $\{2\}$ and so $\alpha'$ is not
contained in any subgroup of ${\cal OP}_{3}$ and nor is $\alpha'$
order-preserving. In the next example $\alpha\in E({\cal O}_{2n})$,
but $\alpha'$ has no fixed points. 

\textbf{Example 3.1.5 }Take $\alpha\in E({\cal O}_{2n})$ so that
$R(\alpha)=F(\alpha)$, putting $R(\alpha)$ as the set of odd members
of $[2n]$ and for each even integer $i\in[2n]$ we put $i\alpha=i+1$.
This yields an order-preserving idempotent $\alpha$ on $[2n]$ of
rank $n$ for which
\[
(0\alpha,1\alpha,\cdots,(2n-2)\alpha,(2n-1)\alpha)=(1,1,3,3,\cdots,2n-1,2n-1).
\]
Now $\alpha=\Phi(K,R,0,1)$ where $K=\{0,2,4,\cdots,2n-2\}$ and \newline $R=\{1,3,5,\cdots,2n-1\}$.
Hence $\alpha'=\Phi(R,K,0,1)$. The kernel classes of $\alpha'$ have
the form $(i,i+1),\,(i=1,3,\cdots,2n-1)$. We see that $\alpha'$
contains the $n$-cycle: 
\[
\sigma=(2n-2\,2n-4\,2n-6\,\cdots\,2\,0),
\]
and for all odd $i$ we have $i\alpha'=i-1$. (The digraph $G(\alpha)$
has exactly one component consisting of the $n$-cycle $\sigma$ along
with $n$ endpoints, one for each point on $\sigma$.) In particular
$c(\alpha)=1$ but $c(\alpha')=n$. 

\subsection{The dual matching involution of ${\cal P}_{n}$}

The mapping on ${\cal P}_{n}$ defined by $\alpha\mapsto\alpha\gamma$
(resp. $\alpha\mapsto\gamma\alpha)$ is an involution on ${\cal P}_{n}$
that maps ${\cal OP}_{n}$ onto ${\cal OR}_{n}$ and maps ${\cal OR}_{n}$
onto ${\cal OP}_{n}$. Additionally, for any $\alpha\gamma\in{\cal OR}_{n}$
(resp. $\gamma\alpha\in{\cal OR}_{n}$) and $\alpha'\in V(\alpha)$
in ${\cal OP}_{n}$, $\gamma\alpha'\in V(\alpha\gamma)$ (resp. $\alpha'\gamma\in V(\gamma\alpha)$)
which lies in ${\cal OR}_{n}$. The upshot of this is that any permutation
matching $(\cdot')$ on ${\cal OP}_{n}$ may be extended to one on
${\cal P}_{n}$ by defining $(\alpha\gamma)'=\gamma\alpha'$ (or dually,
$(\gamma\alpha)'=\alpha'\gamma)$. However, if $(\cdot')$ is an involution
matching, the same is not generally true of either of these extensions,
even in the case of the natural inverse matching. Lemma 2.1.11 supplies
enough information to make this point. 

Let $\beta=\Phi(K,R,i,-1)\in{\cal OR}_{n}$ so that $\rho(\beta)=(K,R)$.
Writing $\beta=\alpha\gamma$ so that $\alpha=\beta\gamma$ we get
$\rho(\alpha)=(K,n-1-R)$ and so $\rho(\alpha')=(n-1-R,K)$. Then
for $\overline{\beta}=\gamma\alpha'$ we have $\rho(\overline{\beta})=(R+1,K)$.
Factorizing $\overline{\beta}$ as $(\gamma\alpha'\gamma)\gamma$
we obtain $\rho(\gamma\alpha'\gamma)=(R+1,n-K-1)$ so that $\rho((\gamma\alpha'\gamma)')=(n-K-1,R+1)$.
Finally, $\overline{\overline{\beta}}=\gamma(\gamma\alpha'\gamma)'$
for which we have $\rho(\overline{\overline{\beta}})=(K+1,R+1)$.
In particular we see that in general $\overline{\overline{\beta}}\neq\beta$,
as $K+1=K$ if and only if $\beta$ is a member of the group of units
of ${\cal P}_{n}$. However, by replacing the standard linear ordering
by the reverse, or as we shall call it the \emph{dual ordering} of
$[n]$, we automatically obtain a dual involution matching on ${\cal P}_{n}$,
which we shall denote by $(\overline{\cdot})$. This generates a distinct
matching involution of ${\cal P}_{n}$ to that of Theorem 3.1.2 and
we now seek to express $(\overline{\cdot})$ in \emph{KRik-}coordinates. 

Let $\alpha=\Phi(K,R,i,k)$. Under $(\overline{\cdot})$, each $r\in R$
is mapped to the initial point of $r\alpha^{-1}$ in the dual ordering,
which is the terminal point $r\alpha^{-1}$ in the standard ordering.
It follows that $X\overline{\alpha}=K-1$. Similarly, under $(\overline{\cdot})$,
$R$ becomes the set of initial points of kernel classes in the dual
ordering, which is then the set of terminal points of those same classes
when expressed in the standard ordering, and so $K(\overline{\alpha})=R+1$.
Therefore $\rho(\overline{\alpha})=(R+1,K-1).$ Since the choice of
ordering does not affect whether or not $\alpha\in{\cal P}_{n}$ preserves
or reverses orientation, we may write $\overline{\alpha}=\Phi(R+1,K-1,i(\overline{\alpha}),k),$
where it only remains to determine $i(\overline{\alpha})$. 

\textbf{Theorem 3.2.2 }The dual matching involution $\overline{(\cdot)}:{\cal P}_{n}\rightarrow{\cal P}_{n}$
acts by $\alpha=\Phi(K,R,i,k)\mapsto\Phi(R+1,K-1),i(\overline{\alpha}),k)$
where: 

Case (1): $0\not\in K$ and $n-1\not\in R$: $i(\overline{\alpha})=1+k(1-i)$; 

Case (2): $0\not\in K$ and $n-1\in R$: $i(\overline{\alpha})=1-ki$;

Case (3): $0\in K$ and $n-1\not\in R$: $i(\overline{\alpha})=k(1-i)$;

Case (4): $0\in K$ and $n-1\in R$: $i(\overline{\alpha})=-ki.$

Moreover, $\alpha$ is in Case (1/4) if and only if $\overline{\alpha}$
is in Case (1/4) and $\alpha$ is in Case (2/3) if and only if $\overline{\alpha}$
is in Case (3/2). 

\textbf{Proof }Since $\rho(\overline{\alpha})=(R+1,K-1)$ it follows
that $\rho(\alpha\overline{\alpha})=(K,K-1)$. Observe that the unique
member of $K-1$ in the kernel class $K_{j}$ is $a_{j+1}-1$ and
since $\alpha\overline{\alpha}\in E({\cal P}_{n})$ we have:
\[
K_{j}\alpha\overline{\alpha}=a_{j+1}-1
\]
\[
\Rightarrow r_{i+kj}\overline{\alpha}=a_{j+1}-1\Rightarrow r_{i+j}\overline{\alpha}=a_{kj+1}-1
\]
\[
\Rightarrow r_{j}\overline{\alpha}=a_{k(j-i)+1}-1=a_{(1-ki)+kj}-1
\]
\begin{equation}
\therefore(r_{j}+1)\overline{\alpha}=r_{j+1}\overline{\alpha}=a_{(1+k(1-i))+kj}-1
\end{equation}
The value of $\overline{i}$ now falls into four cases.

Case (1): $0\not\in K$ and $n-1\not\in R$. Since $1\leq a_{0}$
and $r_{t-1}\le n-2$ we have, in ascending order:
\begin{equation}
K-1=\{a_{0}-1,a_{1}-1,\cdots,a_{t-1}-1\},\,R+1=\{r_{0}+1,r_{1}+1,\cdots,r_{t-1}+1\}
\end{equation}
Then putting $j=0$ in (3) yields $(r_{0}+1)\overline{\alpha}=a_{1+k(1-i)}-1$
and so $\overline{\alpha}=\Phi(R+1,K-1,1+k(1-i),k)$. 

Case (2): $0\not\in K$ but $n-1\in R$ so that the ordered set $R+1=\{r_{t-1}+1=0,r_{0}+1,\cdots,r_{t-2}+1\}$.
Then since $r_{t-1}+1$ is the initial entry of $R+1$ we substitute
$j=t-1\equiv-1$ (mod $t$) into (3) to recover that $\overline{i}$
is $1-ki$ and so $\overline{\alpha}=\Phi(R+1,K-1,1-ki,k)$. 

Case (3): $0\in K$ but $n-1\not\in R$ so that the ordered set $K-1=\{a_{1}-1,a_{2}-1,\cdots,a_{t-1}-1,a_{0}-1=n-1\}$.
As in Case (1) we put $j=0$ in (3) to get $a_{1+k(1-i)}-1$ but since
each entry is now listed one place earlier in the ordered set $K-1$
compared to Case (1), we subtract $1$ from the outcome in Case (1)
to obtain $\overline{\alpha}=\Phi(R+1,K-1,k(1-i),k)$.

Case (4): $0\in K$ and $n-1\in R$ so the ordered set $R+1$ is as
in Case (2) and $K-1$ is as in Case (3). Hence $r_{-1}+1$ is the
first entry of $R+1$, which, by (3), is mapped to $a_{1-ki}$, which
is the entry at position $-ki$ in the list of $K-1.$ Hence we obtain
$\overline{\alpha}=\Phi(R+1,K-1,-ki,k)$. 

Also note that $0\in K\Leftrightarrow n-1\in K-1$ and $n-1\in R\Leftrightarrow0\in R+1$.
It follows that $\alpha$ is in Case (1/4) if and only if $\overline{\alpha}$
is in Case (1/4) and that $\alpha$ is in Case (2/3) if and only if
$\overline{\alpha}$ is in Case (3/2). $\qed$

\textbf{Remarks 3.2.3 }If we take the union of the natural involution
matching $(\cdot)'$ on ${\cal OP}_{n}$ with the dual involution
matching $\overline{(\cdot)}$ on its complement in ${\cal P}_{n}$,
we have another involution matching on ${\cal P}_{n}$. Since the
natural involution matching on ${\cal P}_{n}$ is the unique involution
matching that extends $(\cdot)'$ to ${\cal P}_{n}$ while preserving
${\cal H}$-classes, it follows that this alternative matching is
an example of an involution matching of ${\cal P}_{n}$ that does
not preserve ${\cal H}.$ 

We may check directly that $\overline{(\cdot)}$ defines an involution:
the only non-obvious feature is that the formulas for the shift co-ordinates
are self-inverse in Cases (1) and (4) and mutally inverse for Cases
(2) and (3) but these are readily checked: for example in Case (1):
$1+k(1-(1+k(1-i)))=i.$

An approach by `half duals' leads to permutation matchings that are
however not involutions. For instance we may look to inverses that
map to terminal points of kernel classes while keeping $R$ as the
set of initial points of kernel classes of the inverse. In detail,
for $\alpha\in{\cal P}_{n}$ such that $\rho(\alpha)=(K,R)$ the ${\cal H}$-classes
defined by the kernel-range pairs $(K,K-1)$ and $(R,R)$ are groups
and so there exists a unique inverse $\dot{\alpha}$ of $\alpha$
such that $\rho(\dot{\alpha})=(R,K-1)$. The mapping $\alpha\mapsto\dot{\alpha}$
is then a permutation matching of ${\cal P}_{n}$ but not an involution
as $\rho(\ddot{\alpha})=(K-1,R-1)\neq\rho(\alpha)$ (unless $\alpha$
lies in the group of units of ${\cal P}_{n}$). A dual comment applies
to the other half dual where the inverse of $\alpha$ lies in the
${\cal H}$-class defined by the pair $(R+1,K)$. 

\textbf{Examples 3.2.4 }We take $n=8$, $t=4$, $K=\{0,2,4,6\}$,
$R=\{1,3,5,7\}$ and put $\alpha=\Phi(K,R,0,1)$. We have $R+1=\{0,2,4,6\}=K,\,\,K-1=\{1,3,5,7\}=R.$
Here we have $0\in K$ and $n-1=7\in R$ so that we are in Case (4).
By Theorem 3.2.1 we obtain $\overline{\alpha}=\Phi(R+1,K-1,0,1)=\alpha$,
and indeed $\alpha$ is an idempotent. In contrast, the natural inverse
$\alpha'=\Phi(R,K,0,1)$: 
\[
\alpha=\alpha^{2}=\begin{pmatrix}0 & 1 & 2 & 3 & 4 & 5 & 6 & 7\\
1 & 1 & 3 & 3 & 5 & 5 & 7 & 7
\end{pmatrix}=\overline{\alpha},\,\,\alpha'=\begin{pmatrix}0 & 1 & 2 & 3 & 4 & 5 & 6 & 7\\
6 & 0 & 0 & 2 & 2 & 4 & 4 & 6
\end{pmatrix}.
\]
Next we revisit the first of Examples 3.1.3: $\alpha=\Phi(K,R,3,-1){\cal \in OR}_{8}$,
where $K=\{0,2,4,6\}$\textbf{ }and $R=\{1,3,5,7\}$. Since $0\in K$
and $n-1=7\in R$ we are in Case (4) and so $\overline{\alpha}=\Phi(R+1,K-1,3,-1)$
and so $R+1=\{0,2,4,6\}=K$ and $K-1=\{1,3,5,7\}=R$. Hence: 
\[
\alpha=\begin{pmatrix}0 & 1 & 2 & 3 & 4 & 5 & 6 & 7\\
7 & 7 & 5 & 5 & 3 & 3 & 1 & 1
\end{pmatrix}=\overline{\alpha}.
\]
Since $\overline{\alpha}=\alpha$, therefore $\overline{\overline{\alpha}}=\alpha$
also, and $\alpha^{3}=\alpha.$ In particular, $\overline{\alpha}\neq\alpha'$
. In general, ${\cal OR}_{n}$ contains no idempotents of rank greater
than $2$, so that $\alpha=\alpha^{2}$ is impossible for $\alpha\in{\cal OR}_{n}\setminus{\cal OP}_{n}$. 

As a third example let $n=10$, $t=5$, $K=\{1,3,5,7,8\},$ $R=\{2,3,4,6,8\}$
with $\alpha=\Phi(K,R,4,-1)$. Here $0\not\in K$ and $n-1=9\not\in R$
and so we are in Case (1). Note that since $k=-1$ we have for all
$i$ that $\overline{i}=1-(1-i)=i$, so in particular $i(\overline{\alpha})=4$
and so $\overline{\alpha}=\Phi(R+1,K-1,4,-1)$. We see that $R+1=\{3,4,5,7,9\}$
and $K-1=\{0,2,4,6,7\}$: 
\[
\alpha=\begin{pmatrix}0 & 1 & 2 & 3 & 4 & 5 & 6 & 7 & 8 & 9\\
2 & 8 & 8 & 6 & 6 & 4 & 4 & 3 & 2 & 2
\end{pmatrix}\,\,\overline{\alpha}=\begin{pmatrix}0 & 1 & 2 & 3 & 4 & 5 & 6 & 7 & 8 & 9\\
0 & 0 & 0 & 7 & 6 & 4 & 4 & 2 & 2 & 0
\end{pmatrix}.
\]
Beginning with $\overline{\alpha}=\Phi(R+1,K-1,4,-1)$ we have $0\not\in R+1$
and $9\not\in K-1$ so we are (necessarily) again in Case (1). As
before $i(\overline{\overline{\alpha}})=4$ and we obtain as expected:
\[
\overline{\overline{\alpha}}=\Phi((K-1)+1,(R+1)-1,4,-1)=\Phi(K,R,4,-1)=\alpha.
\]
In contrast the natural inverse of $\alpha$ is $\alpha'=\Phi(R,K,1,-1).$ 

As an example illustrating Cases (2/3) let $n=10,$ $t=6$, $K=\{1,2,5,7,8,9\}$,
$R=\{0,4,5,6,7,9\}$ with $\alpha=\Phi(K,R,4,1)$. Here $0\not\in K$
but $n-1=9\in R$, putting $\alpha$ in Case (2). Theorem 3.2.1 gives
$i(\overline{\alpha})=1-ki(\alpha)=1-4=-3\equiv3$ (mod $6$). Hence
$\overline{\alpha}=\Phi(R+1,K-1,3,1)$, where $0\in R+1=\{0,1,5,6,7,8\}$
and $9\not\in K-1=\{0,1,4,6,7,8\}$, placing $\overline{\alpha}$
in Case (3). Therefore 
\[
\alpha=\begin{pmatrix}0 & 1 & 2 & 3 & 4 & 5 & 6 & 7 & 8 & 9\\
6 & 7 & 9 & 9 & 9 & 0 & 0 & 4 & 5 & 6
\end{pmatrix}\,\,\overline{\alpha}=\begin{pmatrix}0 & 1 & 2 & 3 & 4 & 5 & 6 & 7 & 8 & 9\\
6 & 7 & 7 & 7 & 7 & 8 & 0 & 1 & 4 & 4
\end{pmatrix},
\]
and $\overline{\overline{\alpha}}=\alpha$: $i(\overline{\overline{\alpha}})=k(1-i(\overline{\alpha}))=1(1-3)=-2\equiv4$
(mod $6)=i(\alpha)$. 

\section{Inverse covers and involution matchings for ${\cal T}_{n}$}

\subsection{Inverse covers for ${\cal T}_{n}$ }

In 1971 it was shown by Schein {[}10{]} that every finite full transformation
semigroup ${\cal T}_{n}$ is covered by its inverse subsemigroups,
a result that does not extend to the case of an infinite base set,
{[}5, Ex. 6.2.8{]}. If there existed a cover ${\cal A}=\{A_{i}\}_{1\leq i\leq m}$
of inverse subsemigroups of ${\cal T}_{n}$ with the additional property
that the intersection of any pair of semigroups of ${\cal A}$ was
also an inverse subsemigroup of ${\cal T}_{n}$ then we could deduce
(as explained below) that ${\cal T}_{n}$ had an involution matching.
(The semigroups ${\cal OP}_{n}$ and ${\cal P}_{n}$ of the previous
sections have an inverse cover only if $n\leq3$, {[}2{]}).

It is convenient in what follows to consider the empty set also to
be an inverse semigroup. Suppose that ${\cal A}=\{A_{i}\}_{1\leq i\leq m}$
is an \emph{inverse cover} of ${\cal T}_{n}$ meaning that each $A_{i}$
is an inverse subsemigroup of ${\cal T}_{n}$ and that $\cup_{i=1}^{m}A_{i}={\cal T}_{n}$.
It follows that, for all $1\leq i,j\leq m$, the subsemigroup $A_{i,j}=:A_{i}\cap A_{j}$
of ${\cal T}_{n}$ has commuting idempotents. Indeed it follows easily
from this that Reg$(A_{i,j})$, the set of regular elements of $A_{i,j}$,
forms an inverse subsemigroup of $A_{i,j}$. However it does not automatically
follow that $A_{i,j}=$ Reg$(A_{i,j})$. 

Let $S$ be an arbitrary semigroup and $a\in S$. We say that $b\in V(a)$
is a \emph{strong inverse }of $a$ if the subsemigroup $\langle a,b\rangle$
of $S$ is an inverse semigroup. We denote the set of strong inverses
of $a$ by $S(a)$. We next observe that $S$ has an inverse cover
if and only if every element of $S$ has a strong inverse for, on
the one hand, if every element $a$ has a strong inverse then $S$
is covered by its inverse subsemigroups $\langle a,b\rangle$ where
$b\in S(a)$. On the other hand suppose that $S$ has an inverse cover.
Take $a\in S$ and choose an inverse subsemigroup $A_{a}$ of $S$
containing $a$ and let $b$ be the (unique) inverse of $a$ in $A_{a}$.
Then $A=\langle a,b\rangle$ is a subsemigroup of $A_{a}$ with commuting
idempotents and every element of $A$ is regular as for any product
$p=c_{1}c_{2}\cdots c_{k}\in A$ $(c_{i}\in\{a,b\}$) we see that
$p'=c_{k}'c_{k-1}'\cdots c_{1}'$ is an inverse of $p$ in $A$, where
we take $a'=b$ and $b'=a$, because both products take place within
the inverse semigroup $A_{a}$. It follows that to prove that a given
semigroup $S$ has an inverse cover is equivalent to showing that
$S(a)$ is non-empty for every $a\in S$.

The following general observation applies to any inverse cover ${\cal A}=\{A_{i}\}_{i\in I}$
of an arbitrary semigroup $S$: if the pairwise intersection of any
two members of ${\cal A}$ is an inverse subsemigroup of $S$ then
the same is true of arbitrary intersections. To see this let $J\subseteq I$
and consider $A=\cap_{j\in J}A_{j}$. Either $A$ is the empty inverse
subsemigroup or we may choose $a\in A$ and consider an arbitrary
$A_{j}$ $(j\in J$). Then $a$ has a unique inverse $a^{-1}$ in
$A_{j}$. Now let $k\in J$. By hypothesis, $A_{j}\cap A_{k}$ is
an inverse subsemigroup of $S$ that contains $a$. Since $A_{j}\cap A_{k}$
is an inverse subsemigroup of the inverse semigroup $A_{j}$, it follows
that the unique inverse of $a$ in $A_{j}\cap A_{k}$ is $a^{-1}$.
Since $k\in J$ was arbitrary, it follows that $a^{-1}\in A$ and
so $A$ is indeed an inverse subsemigroup of $S$. We shall say that
$S$ has a \emph{closed inverse cover} if $S$ has a cover by inverse
subsemigroups for which all pairwise intersections of its members
are themselves inverse semigroups.

\textbf{Theorem 4.1.1} For a finite semigroup $S$:

(i) if $S$ has a closed inverse cover then $S$ has an involution
matching by strong inverses.

(ii) If every element $a\in S$ has a unique strong inverse $b$ then
$S$ has a closed inverse cover
\[
{\cal C}=\{\cap_{i=1}^{k}U_{i},\,U_{i}=\langle a,b\rangle,\,a\in S,\,k\geq1\}.
\]

(iii) If $a$ is a group element of $S$ then $a^{-1},$ the group
inverse of $a$ in $S$, is the unique strong inverse of $a$. 

\textbf{Proof }(i) Suppose there exists a closed inverse cover ${\cal A}=\{A_{0}\}_{0\leq i\leq m}$
of $S$ where, without loss, we include $\emptyset$ as $A_{0}$.
The collection ${\cal A}$ is partially ordered by inclusion. Since
every partial order may be extended to a total order, we may order
the members of ${\cal A}$ in such a way that if $A_{i}\subset A_{j}$,
then $A_{i}$ appears before $A_{j}$ in the list. This is assumed
in the following argument.

We now show how ${\cal A}$ could be used to build an involution matching
$(\cdot')$ of $S$ for which $a'\in S(a)$. First $A_{0}$ has an
involution matching $(\cdot')$ in the empty function. Next let $U=\cup_{i=0}^{k}A_{i}$
$(k\geq1)$. Suppose inductively that we have extended the involution
$(\cdot')$ to $V=\cup_{i=0}^{k-1}A_{i}$ and that for each $a\in A_{j}$,
$(0\leq j\leq k-1)$ $a'\in A_{j}$ (so that $a'\in S(a))$. Let $a\in A_{k}$.
Suppose first that $a\in V$ so that $a\in A_{j}$ for some $0\leq j\leq k-1$.
Then $a'$ is already defined and by the nature of the linear order
we have imposed on ${\cal A}$, $A_{i}=A_{j}\cap A_{k}\subseteq V$,
with $i\leq j\leq k-1$. Therefore by induction we have $a'\in A_{j}\cap A_{k}$
and so the induction continues in the case where $a\in A_{k}\cap V$. 

Otherwise $a\not\in V$. Then there exists a unique strong inverse
$a'\in S(a)\cap A_{k}$. What is more $a'\not\in V$ for if to the
contrary $a'\in A_{j}$ say $(0\leq j\leq k-1)$ then $a'$ is again
a member of the inverse semigroup $A_{i}=A_{j}\cap A_{k}\subseteq V$,
where $i\leq k-1$. In this event, $(a')'$ is already defined and
would be an inverse of $a'$ in $A_{j}\cap A_{k}$, whence $(a')'=a$.
But then $a\in A_{j}$, contrary to our choice of $a$. It follows
that $a'\not\in V$ and so we may extend the involution by strong
inverses $(\cdot')$ to $U=V\cup A_{k}$ by setting $a'$ as the unique
inverse of $a$ in $A_{k}\setminus V$ for all $a\in A_{k}\setminus V$.
Therefore we see that in both cases the induction continues. Since
${\cal A}$ covers $S$, the process terminates when $k=m$, yielding
an involution matching by strong inverses $(\cdot')$ of $S$. 

(ii) Let $U,V\in{\cal A}$ and suppose that $U\cap V\neq\emptyset$.
For any $a\in U\cap V$ let $b$ be the unique strong inverse of $a$.
Let $u\in V(a)$ in the inverse semigroup $U$. Then $u\in S(a)$,
whence $u=b$. We may draw the corresponding conclusion for the inverse
$v\in V(a)\cap V$, so that $u=b=v$. In particular $b\in U\cap V$,
whence it follows that $U\cap V$ is an inverse subsemigroup. Therefore
by adjoining all intersections $U_{1}\cap U_{2}\cap\cdots\cap U_{k}$
$(k\geq2)$ of members $U_{i}\in{\cal A}$ to the inverse cover ${\cal A}$
we generate a closed inverse cover for $S$. 

(iii) Clearly $a^{-1}\in S(a)$. Consider an arbitrary $b\in S(a)$
and let $e=ab,$ $f=ba$. Then we have $e{\cal R}a{\cal L}f$ in $S$.
Then since $H_{a}$ is a group we have $fe{\cal H}b$. Since $b\in S(a)$
it follows that $ef=fe\in E(S)$. But then $a{\cal H}ef=fe{\cal H}b$
and so $b,a^{-1}\in V(a)$ with $b{\cal H}a^{-1}$, whence $b=a^{-1}$
by uniqueness of inverses within an ${\cal H}$-class. $\qed$

\textbf{Remark 4.1.2} As explained prior to Theorem 4.1.1, for any
semigroup with a closed inverse cover, the intersection of any collection
of members of ${\cal A}$ is also an inverse subsemigroup of $S$
and so we may assume further that ${\cal A}$ is closed under the
taking of arbitrary intersections of its members. This allows the
argument of the previous proof to be extended to arbitrary semigroups
through the Axiom of choice and transfinite induction. Part (ii) of
Theorem 4.1.1 is a partial converse of part (i). It remains open as
to whether or not the full converse holds. 

In the next section, we shall prove that in general ${\cal T}_{n}$
has no involution matching by strong inverses, from which it follows
from the contrapositive of Theorem 4.1.1(i) that ${\cal T}_{n}$ has
no closed inverse cover. 

\newpage

\subsection{Closed inverse covers for ${\cal T}_{n}$ do not exist}

In this section our context throughout will be ${\cal T}_{n}$. The
account here of the construction of strong inverses in ${\cal T}_{n}$
follows {[}5, Section 6.2{]}. 

Let $\alpha\in{\cal T}_{n}$. For $x\in X_{n}$ the \emph{depth }of
$x$, denoted by $d(x)$, is the length of the longest dipath in $G(\alpha)$
ending at $x$; if $x\in$ stran$(\alpha)$ we conventionally define
$d(x)=\infty$. Note that $d(x)=k<\infty$ if and only if $x\in X\alpha^{k}\setminus X\alpha^{k+1}$.
The \emph{height }of $x$, denoted by $h(x)$ is the least positive
integer $k$ such that $d(x\alpha^{k})\geq d(x)+k+1$; again we take
$h(x)=\infty$ if $x\in$ stran$(\alpha)$. The height of $x$ is
the length of the dipath which begins at $x$ and terminates at the
first point $u$ which is also the terminal point of some dipath that
is longer than the dipath from $x$ to $u$. A necessary condition
for membership of $S(a)$ is the following.

\textbf{Lemma 4.2.1} Let $\beta\in S(\alpha)$. Then for all $x\in X\alpha$,
$x\beta$ is a member of $x\alpha^{-1}$ of maximal depth.

When constructing strong inverses, the correct treatment of the endpoints
of $G(\alpha),$ which are those $x\in X_{n}$ for which $d(x)=0$,
is more complicated. The next parameter is defined on the vertices
of $G(\alpha)$ in terms of some fixed but arbitrary $\beta\in{\cal T}_{n}$,
but is only significant when $\beta\in V(\alpha)$. For each $x\in X_{n}$
the \emph{grasp }$g(x)$ of $\alpha$ is the greatest non-negative
integer $k$ such that $x\alpha^{k}\beta^{k}=x$. 

\textbf{Lemma 4.2.2 }Let $\beta\in S(\alpha)$. If $d(x)=0$, and
$h(x)=h$ then $x\beta=y$ satisfies $g(y)\geq g(x)+1$ and $y\alpha^{h+1}=x\alpha^{h}$. 

Lemmas 4.2.1 and 4.2.2 are all we require here. However, if $\beta\in{\cal T}_{n}$
satisfies these conditions together with the equality $x\beta\alpha\cdot\alpha^{g(x)+1}\beta^{g(x)+1}=x\alpha^{g(x)+1}\beta^{g(x)+1}\cdot\beta\alpha$,
it may then be proved that that $\beta\in S(\alpha)$. We may show
from this point that ${\cal T}_{n}$ has an inverse cover for the
lemmas represent the two stages in the construction of a particular
type of strong inverse $\beta\in S(\alpha)$: Lemma 4.2.1 applies
to points of positive depth in $G(\alpha)$, while for each endpoint
$x$ we may follow the dipath (of length $k$ say) from $x$ until
we meet a point $u$ of depth exceeding $k$. Then $u\beta^{t}$ has
already been defined for all $0\leq t\leq k+1$ and we then put $x\beta=x\alpha^{k}\beta^{k+1}$.
This stage can always be carried out and indeed this $\beta\in S(\alpha)$
is uniquely determined by the choices made in determining $X\alpha\beta$.
The outcome of this is a particular strong inverse $\beta\in S(\alpha)$
for such a $\beta$ will also satisfy the additional condition and
indeed the set of all idempotents $\alpha^{t}\beta^{t},$$\beta^{s}\alpha^{s}$
then commute with each other, as is required for $\langle\alpha,\beta\rangle$
to be an inverse semigroup. 

The main result of this section is the following.

\textbf{Theorem 4.2.3 }The full transformation semigroup ${\cal T}_{n}$
has a closed inverse cover if and only if $n\leq3$. 

\textbf{Lemma 4.2.4 }For $n\leq3$, ${\cal T}_{n}$ has a closed inverse
cover. 

\textbf{Proof }For $n=1,2$ we note that ${\cal T}_{n}$ is a union
of groups so the claim follows from Theorem 4.1.1. Although ${\cal T}_{3}$
is not a union of groups, we may verify that each $\alpha\in{\cal T}_{3}$
has a unique strong inverse as follows.

In general, an element $\alpha\in{\cal T}_{n}$ is a group element
if and only if $X\alpha=X\alpha^{2}$. It follows that all members
of ${\cal T}_{3}$ of ranks 1 or 3 are group elements. There are $3^{3}-3-3!=18$
mappings in ${\cal T}_{3}$ of rank $2$. All of those with two components
are idempotent (these number $3\times2=6$). Those with one component
for which $|X\alpha^{2}|=2$ are group elements (these also number
$6$). This leaves $18-6-6=6$ mappings of rank $2$ with a single
component and for which $|X\alpha^{2}|=1$. These are evidently the
$6$ mappings $\alpha$ of the form $a\mapsto b\mapsto c\mapsto c$
where $\{a,b,c\}=\{1,2,3\}$, which we denote for our current purposes
by $(a\,b\,c)$. Observe that each such $\alpha$ has a unique strong
inverse, which is $\alpha'=(b\,a\,c)$. The result now follows by
Theorem 4.1.1(ii). $\qed$ 

\textbf{Examples 4.2.5 }For $\alpha=(1\,2\,3)$ we have $S(\alpha)=\alpha'=(2\,1\,3)$,
$\alpha\alpha':\text{\ensuremath{1\mapsto1,\,\,2,3\mapsto3}}$, $\alpha'\alpha:1,3\mapsto3,\,2\mapsto2$.
The inverse subsemigroup $U_{3}=\langle\alpha,\alpha'\rangle$ is
a 5-element combinatorial Brandt semigroup with zero element given
by $\alpha^{2}=\alpha'^{2}$, which is the constant mapping with range
$\{3\}$. The subsemigroups $U_{1},U_{2},U_{3}$ are pairwise disjoint.
However not all intersections of distinct members of the inverse cover
${\cal C=}\{\langle a,b\rangle:b\in S(a)\}$ are empty: for example
consider the mapping $\gamma:1\mapsto3,\,2,3\mapsto1$, which is its
own strong inverse. Since $\gamma^{2}=\alpha\alpha'$ we obtain $\langle\gamma\rangle\cap\langle\alpha,\alpha'\rangle=\{\alpha\alpha'\}$. 

By way of contrast, let us examine a subsemigroup $\langle\alpha,\alpha_{1}\rangle$,
where $\alpha_{1}\in V(\alpha)\setminus S(\alpha)$: we take $\alpha_{1}=(3\,2\,1)$.
Then $e=\alpha\alpha_{1}\in E({\cal T}_{3})$ and has fixed point
set of $\{1,2\}$ with $3e=2$. Also $f=\alpha^{2}\in E({\cal T}_{3})$
is the constant mapping with range $\{3\}$, so that $ef=f$. However
$fe$ is the constant mapping onto $2$ and so idempotents do not
commute in $U=\langle\alpha,\alpha_{1}\rangle$. In fact $U$ is a
7-element regular subsemigroup of $S$, containing all three constant
mappings, as $\alpha_{1}^{2}$ is the constant with range $\{1\}$. 

Given Lemma 4.2.4 we now need to prove that ${\cal T}_{n}$ does not
have a closed inverse cover for $n\geq4$. The remaining substantial
task is to show that ${\cal T}_{4}$ has no involution matching through
strong inverses for that implies that ${\cal T}_{4}$ has no closed
inverse cover. For $n\geq5$ we then consider the copy of ${\cal T}_{4}$
embedded in ${\cal T}_{n}$ defined by $T=\{\alpha\in{\cal T}_{n}:\,k\alpha=k\,\forall\,k\geq5\}$.
Suppose that ${\cal C}$ were a closed inverse cover for ${\cal T}_{n}$.
If $\alpha\in T$ then for any $\beta\in S(\alpha)$ we have $\beta\in T$.
It follows that ${\cal C}_{T}=\{A\cap T:A\in{\cal C}\}$ would be
a closed inverse cover for $T$, which is isomorphic to ${\cal T}_{4}$,
which would then also have such a cover. Therefore, to complete the
proof of our theorem, it remains only to show that ${\cal T}_{4}$
does not have an involution matching by strong inverses. 

First we identify every member of ${\cal T}_{4}$ that possesses a
unique strong inverse, a collection that includes all mappings of
ranks 1 or $4$ as these are group elements. Indeed for any rank we
only need consider non-group elements, which are the mappings $\alpha$
such that $X\alpha^{2}$ is proper subset of $X\alpha$. 

Consider mappings of rank $3$. It follows that $|X\alpha^{2}|\leq2$.
If $\alpha$ has two components, since $X\alpha^{2}\neq X\alpha,$
it follows that $\alpha$ has an isolated fixed point $d$ say and
a second component of the form $a\mapsto b\mapsto c$, which has a
unique strong inverse $b\mapsto a\mapsto c$, $d\mapsto d$. If $\alpha$
has just one component with $|$stran$(\alpha)|=2$ then $\alpha$
necessarily now has the form $a\mapsto b\mapsto c\mapsto d\mapsto c$,
which has a unique strong inverse, which is the mapping $b\mapsto a\mapsto d\mapsto c\mapsto d$.
We conclude that all mappings of ranks 1, 3, or 4 have a unique strong
inverse. 

Finally consider mappings of rank $2$. Since we may assume there
exists a point $x\in X\alpha\setminus$stran$(\alpha)$ (as otherwise
$\alpha$ is a group element) it follows that we are restricted to
mappings $\alpha$ with a single component and that component has
a fixed point. The two remaining cases are:

A: the form of a `Y': $\alpha=\begin{pmatrix}a & b & c & d\\
c & c & d & d
\end{pmatrix}$ or B: the form $\beta=\begin{pmatrix}a & b & c & d\\
b & d & d & d
\end{pmatrix}$. 

We next act the strong inverse operator $S(\cdot)$ for mappings of
these two types. We will see that this generates a set of four $9$-cycles,
with each mapping within a cycle sharing the same fixed point. Consequently
these cycles are pairwise disjoint. 

The given mapping $\alpha$ of type A has exactly two strong inverses,
both of which are of type B:
\begin{equation}
\beta_{1}:\,\begin{pmatrix}a & b & c & d\\
d & d & a & d
\end{pmatrix}\,\beta_{2}=\begin{pmatrix}a & b & c & d\\
d & d & b & d
\end{pmatrix}\,(B)
\end{equation}
The mapping $\beta$ of type B also has exactly two strong inverses,
the first of type A, the second of type B:
\begin{equation}
\alpha_{1}:\,\begin{pmatrix}a & b & c & d\\
d & a & a & d
\end{pmatrix}\,(A)\,\,\beta{}_{2}:\begin{pmatrix}a & b & c & d\\
d & a & d & d
\end{pmatrix}\,(B)
\end{equation}
Consider the collection $C$ of all mappings of rank $2$ with two
strong inverses and a common fixed point, $d$. There are $3$ mappings
of type A and $6$ of type B, so that $C$ has $9$ members. We use
the symbols $\alpha$ and $\beta$, with appropriate subscripts, to
denote mappings of types A and of B respectively. 

The strong inverse operator $S(\cdot)$ acting on a point in $C$
outputs exactly two distinct mappings, which are also members of $C$,
in accord with the rules (5) and (6). Let us write $\alpha_{1}=\alpha$
for the type A mapping above. We write $\beta_{1}\rightarrow\alpha_{1}\rightarrow\beta_{2}$
with the arrow indicating the first map is a strong inverse of the
second (so that the reverse arrow is equally valid). We now act the
operator $S(\cdot)$ on the rightmost member of our sequence, which
will produce as outputs the previous member and a new sequence member.
Bearing in mind rules (5) and (6) our sequence $C$ will thus take
on the form:
\begin{equation}
C:\,\beta_{1}\rightarrow\alpha_{1}\rightarrow\beta_{2}\rightarrow\beta_{3}\rightarrow\alpha_{2}\rightarrow\beta_{4}\rightarrow\beta_{5}\rightarrow\alpha_{3}\rightarrow\beta_{6}\cdots.
\end{equation}
The output of $S(\gamma)$ when acting on $\gamma\in C$ comprises
two distinct mappings, neither of which is $\gamma$, and one of which
is the predecessor of $\gamma$ in the sequence. Eventually the output
$S(\gamma)$ will produce a repeated member of $C$ (in addition to
the predecessor of $\gamma$), which must appear at least two steps
before $\gamma$. However, all such members of $C$, apart from $\beta_{1}$,
have already had their two strong inverses appear in $C$, and so
cannot have $\gamma$ as a third strong inverse. Therefore the repeated
sequence member is necessarily $\beta_{1}$. Hence $C$ is a cycle
of length $l$ say with $l\geq2$ and $l|9$, and so $l=3$ or $l=9$.
However $l=3$ would imply that $\beta_{1}$ and $\beta_{2}$ were
mutual inverses, which is not the case. Therefore $l=9$ and and the
cycle $C$ is completed by $\beta_{6}\rightarrow\beta_{1}$. 

Suppose now that ${\cal T}_{4}$ possessed an involution $(\cdot')$
by strong inverses. Any mapping $\alpha$ with a unique strong inverse
$\beta$ is necessarily paired with $\beta$ under $'$. This includes
all mappings in ${\cal T}_{4}$ except for the mappings which are
the vertices of the four disjoint $9$-cycles we have just identified.
Each member of such a $9$-cycle $C$ is then paired with an adjacent
partner in that cycle, but since $9$ is odd, this is not possible
and so we have a contradiction. Therefore ${\cal T}_{4}$ has no involution
by strong inverses, which implies by Theorem 4.1.1(i) that ${\cal T}_{4}$
has no inverse cover closed under the taking of intersections. $\qed$

\textbf{Remarks 4.2.6 }We may explicitly calculate the $9$-cycle
$C$ that contains the mapping $\alpha$ above, denoted here as $\alpha_{1}$,
through repeated use of rules (5) and (6) as follows. We write $S(\alpha_{1})=\{\beta_{1},\beta_{2}\}$,
with the $\beta_{i}$ as given in (5). Then following (7) the subsequent
members of $C$ are $\beta_{2}\rightarrow\beta_{3}=\begin{pmatrix}a & b & c & d\\
d & c & d & d
\end{pmatrix}\rightarrow$ $\alpha_{2}=\begin{pmatrix}a & b & c & d\\
b & d & b & d
\end{pmatrix}\rightarrow$ $\beta_{4}=\begin{pmatrix}a & b & c & d\\
d & a & d & d
\end{pmatrix}\rightarrow$ $\beta_{5}=\begin{pmatrix}a & b & c & d\\
b & d & d & d
\end{pmatrix}\rightarrow$$\alpha_{3}=\begin{pmatrix}a & b & c & d\\
d & a & a & d
\end{pmatrix}\rightarrow$ $\beta_{6}=\begin{pmatrix}a & b & c & d\\
c & d & d & d
\end{pmatrix}\rightarrow$ $\beta_{7}=\beta_{1}=\begin{pmatrix}a & b & c & d\\
d & d & a & d
\end{pmatrix}$, giving the anticipated $9$-cycle $C$. 

Running down the ranks from $4$ to $1$, elementary combinatorial
considerations give that:
\[
|E({\cal T}_{4})|=1+2\binom{4}{2}+\Big(3\binom{4}{1}+(2)(2)\binom{4}{2}\Big)+4=1+12+36+4=53.
\]
In a similar fashion, bracketing term sum contributions from a common
rank, the number of non-idempotent self-inverse elements is given
by:
\[
\Big(\binom{4}{2}+\frac{1}{2}\binom{4}{2}\Big)+(2)(3)\binom{4}{2}+2\binom{4}{2}=9+36+12=57.
\]
The number of mappings with a distinct unique strong inverse is given
by:
\[
(3!+2!\binom{4}{3}\big)+\Big((3)(2!)\binom{4}{3}+(2)(2)\binom{4}{2})+4!)+(2)(2)\binom{4}{2}\Big)=14+72+24=110.
\]
The number of mappings with exactly two strong inverses is $4\times9=36$,
giving the total of $(53+57)+110+36=256=4^{4}=|{\cal T}_{4}|$. The
graph of strong inverses of ${\cal T}_{4}$ then consists of $110$
singletons, $55$ pairs, and four $9$-cycles. In particular the above
analysis shows that ${\cal T}_{4}$ does possess a permutation matching
by strong inverses. We may use one of these permutations to construct
an involution matching of ${\cal T}_{4}$. (There are $2^{4}=16$
such permutations, determined by the $2$ choices of orientation of
the $4$ cycles). First consider the $9$-cycle explicitly calculated
above in which all mappings fix a point $d$. The mapping $\alpha=\alpha_{1}$
has an idempotent inverse $\varepsilon_{d}=\begin{pmatrix}a & b & c & d\\
b & b & b & d
\end{pmatrix}$. We then remove the pair $(\varepsilon_{d},\varepsilon_{d})$ from
our permutation, replacing it by $(\varepsilon_{d},\alpha)$ and pair
up the remaining $8$ members of the associated $9$-cycle in neighbouring
pairs. We repeat this procedure with the other three cycles, noting
that there is no repetition of idempotents used in our pairings. This
then yields an involution matching for ${\cal T}_{4}$. 

We close with an example showing however that in general ${\cal T}_{n}$
does not possess a permutation matching by strong inverses. 

\textbf{Example 4.2.7 }Consider the following pair of members of ${\cal T}_{8}$:
\[
\alpha_{1}=\begin{pmatrix}1 & 2 & 3 & 4 & 5 & 6 & 7 & 8\\
2 & 3 & 4 & 5 & 5 & 3 & 8 & 4
\end{pmatrix}\,\,\alpha_{2}=\begin{pmatrix}1 & 2 & 3 & 4 & 5 & 6 & 7 & 8\\
2 & 3 & 4 & 5 & 6 & 8 & 3 & 4
\end{pmatrix}.
\]
The two mappings are identical except for the interchange of the images
of $6$ and $7$, and so their digraphs are isomorphic. They share
a common range: $X\alpha_{1}=X\alpha_{2}=\{2,3,4,5,8\}$. Moreover
for each $x\in X\alpha_{i}$ $(i=1,2)$ there is a unique member of
$y\in x\alpha_{i}^{-1}$ of maximal depth and so by Lemma 4.2.1 we
see that any strong inverse $\beta_{i}\in S(\alpha_{i})$ has the
following form: 
\[
\beta_{1}=\begin{pmatrix}1 & 2 & 3 & 4 & 5 & 6 & 7 & 8\\
- & 1 & 2 & 3 & 5 & - & - & 7
\end{pmatrix}\,\,\beta_{2}=\begin{pmatrix}1 & 2 & 3 & 4 & 5 & 6 & 7 & 8\\
- & 1 & 2 & 3 & 5 & - & - & 6
\end{pmatrix}.
\]
In each case the points of depth zero are $1,6,$ and $7$. For both
mappings and for any strong inverses $\beta_{i}$ we see that $g(1)=3$
so that for any choice of $y=1\beta_{i}$ we have by Lemma 4.2.2 that
$g(y)\geq4$, which implies that $y=1\beta_{i}=5$. To determine $6\beta_{1}$
we note that $6\alpha_{1}=3$ and so $d(6\alpha_{1})=2>0+1=d(6)+1$;
hence $h(6)=1$ and $g(6)=0$. Writing $y=6\beta_{1}$ we have by
Lemma 2.2 that $g(y)\geq1$ and $y\alpha_{1}^{2}=6\alpha_{1}=3$ so
that $y=6\beta_{1}=1$. By the same argument with $6$ replaced by
$7$ we obtain $7\beta_{2}=1$. 

Finally consider $7\beta_{1}$. We have $7\alpha_{1}=8$ so we see
that $d(7\alpha_{1})=1$ and $g(7)=1$, $h(7)=2$ as $d(7\alpha_{1}^{2})=d(4)=3>2+0=2+d(7)$
while $d(7\alpha)=d(8)=1\not>1+0$. Hence we have by Lemma 4.2.2 that
$y=7\beta_{1}$ must satisfy $g(y)\geq2$ and $y\alpha^{3}=7\alpha^{2}=4$
so that $7\beta_{1}=1$. By symmetry we also obtain $6\beta_{2}=1$.
Therefore each of the $\alpha_{i}$ has a unique strong inverse $\beta_{i}$:
\[
\beta_{1}=\begin{pmatrix}1 & 2 & 3 & 4 & 5 & 6 & 7 & 8\\
5 & 1 & 2 & 3 & 5 & 1 & 1 & 7
\end{pmatrix}\,\,\beta_{2}=\begin{pmatrix}1 & 2 & 3 & 4 & 5 & 6 & 7 & 8\\
5 & 1 & 2 & 3 & 5 & 1 & 1 & 6
\end{pmatrix}.
\]
We now consider a third mapping $\beta\in{\cal T}_{8}$ and a putative
strong inverse $\beta'\in S(\beta)$. As before we have $X\beta=\{2,3,4,5,8\}$
and again Lemma 4.2.2 gives the following unique partial definition
of $\beta'$: 
\[
\beta=\begin{pmatrix}1 & 2 & 3 & 4 & 5 & 6 & 7 & 8\\
2 & 3 & 4 & 5 & 5 & 8 & 8 & 4
\end{pmatrix}\,\,\beta'=\begin{pmatrix}1 & 2 & 3 & 4 & 5 & 6 & 7 & 8\\
- & 1 & 2 & 3 & 5 & - & - & -
\end{pmatrix}.
\]
We see that $8\beta'\in\{6,7\}$; for the moment let us make the choice
$8\beta'=7$ and henceforth denote $\beta'$ by $\beta_{1}'$. The
points of zero depth are again $1,6,$ and 7 and the same analysis
that applied to the $\alpha_{i}$ again yields $1\beta_{1}'=5$. Next
we note that $g(6)=0$ and $h(6)=2$ as $d(6\alpha^{2})=d(4)=3>2+0$
but $d(6\alpha)=d(8)=1\not>1+0$. Hence $y=6\beta_{1}'$ must satisfy
$g(y)\geq1$ and $y\beta^{3}=6\beta^{2}=4$ so that $y=6\beta_{1}'=1$.
Finally we have $g(7)=1$ and $h(7)=2$ as for $h(6)$. Hence $y=7\beta'$
must satisfy $g(y)\geq2$ and $y\beta^{3}=6\beta^{2}=4$ so that $y=7\beta'=1$
also. We have then identified one strong inverse of $\beta_{1}'\in S(\beta)$.
Similarly there exists a second strong inverse $\beta_{2}'\in S(\beta)$
determined by the alternative choice $8\beta_{2}'=6$. Exchanging
the roles of the symbols $6$ and $7$ makes no difference to the
images of the other domain points in that we again obtain that $1\beta_{2}'=5,$
$6\beta_{2}'=7\beta_{2}'=1$. Therefore we find that
\[
\beta'_{1}=\beta_{1}=\begin{pmatrix}1 & 2 & 3 & 4 & 5 & 6 & 7 & 8\\
5 & 1 & 2 & 3 & 5 & 1 & 1 & 7
\end{pmatrix}\,\,\beta'_{2}=\beta_{2}=\begin{pmatrix}1 & 2 & 3 & 4 & 5 & 6 & 7 & 8\\
5 & 1 & 2 & 3 & 5 & 1 & 1 & 6
\end{pmatrix}.
\]
The upshot of all this is that we have a set of three members of ${\cal T}_{8}$
in $U=\{\alpha_{1},\alpha_{2},\beta\}$ such that the set $S(U)$
of all strong inverses of elements of $U$ is the two-element set
$S(U)=\{\beta_{1},\beta_{2}\}$. It follows that there is no permutation
matching $(\cdot)'$ on ${\cal T}_{8}$ that maps the set $U$ into
the set $S(U)$, thereby giving us the result mentioned earlier, which
we now formally state.

\textbf{Corollary 4.2.8 }There is no permutation matching $(\cdot')$
of ${\cal T}_{n}$ $(n\geq8)$ such that $a'$ is a strong inverse
of $a$ for all $a\in{\cal T}_{n}$.

\textbf{Proof }Example 4.2.7 shows the corollary is true for $n=8$.
For $n\geq9$ we may extend the above example with each of the the
mappings $\alpha_{1},\alpha_{2},\beta$ acting identically on all
integers exceeding $8$. Since any strong inverse preserves components
we again obtain the conclusion that $S(\alpha_{1})=\{\beta_{1}\}$,
$S(\alpha_{2})=\{\beta_{2}\}$ and $S(\beta)=\{\beta_{1},\beta_{2}\}$,
which implies the result. $\qed$

\end{document}